\documentclass[12pt]{article}
\usepackage{amsfonts,amsmath}
\usepackage{epsfig,graphics}

\renewcommand{\appendix}{%
\renewcommand{\section}{%
\newpage
\thispagestyle{plain}%
\secdef\Appendix\sAppendix}%
\setcounter{section}{0}%
\renewcommand{\thesection}{\Alph{section}}%
}
\newcommand{\Appendix}[2][?]{%
\refstepcounter{section}%
\addcontentsline{toc}{Addendum}%
{\protect\numberline{\appendixname~\thesection}#1}%
{\flushleft\LARGE\bfseries\appendixname\ \thesection\par
\centering#2\par}%
\sectionmark{#1}\vspace{\baselineskip}}

\newcommand{\sAppendix}[1]{%
{\flushright\large\bfseries\appendixname\par
\centering#1\par}%
\vspace{\baselineskip}}

\def\be{\begin{equation}}

\def\bea{\begin{eqnarray}}

\def\eea{\end{eqnarray}}

\renewcommand{\theequation}{\thesection.\arabic{equation}}

\normalbaselines \oddsidemargin 0.5cm \evensidemargin 0.5cm
\begin{document}

\pagestyle{empty}

\begin{center}

\textsf{\Huge {\bf Higher Dimensional Unitary Braid Matrices:
Construction, Associated Structures and Entanglements}}

\vspace{9mm}

{\bf \Large B. Abdesselam$^{a,}$\footnote{Email:
boucif@cpht.polytechnique.fr and  boucif@yahoo.fr}, A.
Chakrabarti$^{b,}$\footnote{Email:
chakra@cpht.polytechnique.fr},\\
V.K. Dobrev$^{c,d,}$\footnote{Email: dobrev@inrne.bas.bg} and S.G.
Mihov$^{c,}$\footnote{Email: smikhov@inrne.bas.bg}}

\vspace{5mm}

  \emph{$^a$ Laboratoire de Physique Quantique de la
Mati\`ere et de Mod\'elisations Math\'ematiques, Centre
Universitaire de Mascara, 29000-Mascara, Alg\'erie\\
and \\
Laboratoire de Physique Th\'eorique, Universit\'e d'Oran
Es-S\'enia, 31100-Oran, Alg\'erie}
  \\
  \vspace{3mm}
  \emph{$^b$ Centre de Physique Th{\'e}orique, Ecole Polytechnique, 91128 Palaiseau Cedex, France.}
  \\
  \vspace{3mm}
  \emph{$^c$ Institute of Nuclear Research and Nuclear Energy
Bulgarian Academy of Sciences 72 Tsarigradsko Chaussee, 1784
Sofia, Bulgaria}\\
 \vspace{3mm}
 \emph{$^d$ Abdus Salam International Center for Theoretical Physics
Strada Costiera 11, 34100 Trieste, Italy}

\end{center}

\begin{abstract}
{\small \noindent We construct $(2n)^2\times (2n)^2$ unitary braid
matrices $\widehat{R}$ for $n\geq 2$ generalizing the class known
for $n=1$. A set of $(2n)\times (2n)$ matrices $\left(I,J,K,L
\right)$ are defined. $\widehat{R}$ is expressed in terms of their
tensor products (such as $K\otimes J$), leading to a canonical
formulation for all $n$. Complex projectors $P_{\pm}$ provide a
basis for our real, unitary $\widehat{R}$. Baxterization is
obtained. Diagonalizations and block-diagonalizations are
presented. The loss of braid property when $\widehat{R}$ $(n>1)$
is block-diagonalized in terms of $\widehat{R}$ $(n=1)$ is pointed
out and explained. For odd dimension $(2n+1)^2\times (2n+1)^2$, a
previously constructed braid matrix is complexified to obtain
unitarity. $\widehat{R}\mathrm{LL}$- and $\widehat{R}
\mathrm{TT}$-algebras, chain Hamiltonians, potentials for
factorizable $S$-matrices, complex non-commutative spaces are all
studied briefly in the context of our unitary braid matrices.
Turaev construction of link invariants is formulated for our case.
We conclude with comments concerning entanglements. }
\end{abstract}

\pagestyle{plain} \setcounter{page}{1}

%\tableofcontents

\newpage

\section{Introduction}
\setcounter{equation}{0}

The $4\times 4$ unitary braid matrix
\begin{equation}
\widehat{R}=\frac 1{\sqrt{2}}\left|\begin{array}{cccc}
   1 & 0 & 0 & 1 \\
   0 & 1 & -1 & 0\\
   0 & 1 & 1 & 0 \\
   -1 & 0 & 0 & 1 \\
\end{array}\right|
\end{equation}
and the associated $\mathrm{SO3}$ algebra have been studied
extensively in our previous papers (refs. \cite{R1,R2} provide
further sources). Recently this matrix (along with its conjugate
by the permutation matrix $P$) has been studied as the Bell-matrix
in the context of quantum entanglements (see \cite{R3,R4,R5} and
the references therein). Algebraic aspects have been studied
\cite{R4,R5} using different approaches. Joint presence of quantum
and topological entanglements is the theme of ref. \cite{R3}. The
$q$-deformation of such matrix with $q$ at root of unity or
generic and its relationship with quantum computing was also
discussed in refs. \cite{SB,ZKG}. Here we present direct
generalizations of (1.1) to higher dimensions, namely {\it
unitary} matrices of dimensions $(2n)^2\times (2n)^2$ ($n\geq 1$),
which satisfy the braid equation
\begin{equation}
\widehat{R}_{12}\widehat{R}_{23}\widehat{R}_{12}=\widehat{R}_{23}\widehat{R}_{12}\widehat{R}_{23},
\end{equation}
where  $\widehat{R}_{12}=\widehat{R}\otimes I$ and
$\widehat{R}_{23}=I\otimes \widehat{R}$. For $n=1$ one obtains
(1.1). Entanglements are commented upon in sec. 11, referring to
previous sections\footnote{Notations: Our notations will be
systematically defined in the following sections. But here we
point out that our $\widehat{R}$, the braid matrix is the
Yang-Baxter (YB)-matrix of ref. \cite{R3} and our YB-matrix
$R=P\widehat{R}$ is called {\it algebraic YB-matrix} in ref.
\cite{R3}. Other references cited also use different notations.
This should be noted to avoid confusion.}.

Different classes of higher dimensional braid matrices were
constructed and studied in a series of previous papers
\cite{R6,R7,R8,R9}. But unitary was not sought before, though for
the Baxterized form the constraint
\begin{equation}
\widehat{R}\left(\theta\right)\widehat{R}\left(-\theta\right)=I
\end{equation}
is sometimes labeled {\it unitarity}. In the constructions to
follow one has always
\begin{equation}
\widehat{R}^+\widehat{R}=I_{(2n)^2}.
\end{equation}
Moreover, the spectral parameter $\theta$ (or $z=\tanh\theta$) is
introduced in such a way when Baxterizing that
\begin{equation}
\widehat{R}^+\left(\theta\right)=\widehat{R}\left(-\theta\right)
\end{equation}
and unitarity, now really coinciding with (1.3), is maintained.

Our construction is presented in a canonical form for all $n$ by
introducing a set of $2n\times 2n$ matrices $(I,J,K,L)$ with
particularly simple properties and implementing their tensor
products ($K\otimes J$ and so on). The verification of the braid
equation and its Baxterized form now become transparent (see
(2.11), (2.16) and (2.17)). Denoting the $\left(2n\right)^2\times
\left(2n\right)^2$ braid matrix as $\widehat{R}_{(2n)}$ (so that
(1.1) is now $\widehat{R}_{(2)}$) we verify explicitly the
non-equivalence of $\widehat{R}_{(2n)}$ with a block-diagonal,
direct sum of $\widehat{R}_{(2)}$'s in the following, precise
sense: On can indeed construct a matrix, say, $V$ such that
$V\widehat{R}_{(2n)}V^{-1}$ is block-diagonalized into
$\widehat{R}_{(2)}$ blocks. But $V$ does not have a tensored
structure of the type $W\otimes W$ and does not conserve the braid
property (1.2) (see (4.9-12)). The tensored structure displays
this negative result also with clarity. In section 5 we point out
briefly that a certain type of complexification of the odd
dimensional $\widehat{R}\left(\theta\right)$'s of ref. \cite{R6}
leads to unitarity, where however the Baxterized form is
essential. Elsewhere, in this paper, only even dimensional
$\widehat{R}$ is considered. After presenting our general
constructions we study various aspects of our braid matrices in
successive sections: Baxterization, $\widehat{R}\mathrm{LL}$- and
$\widehat{R}\mathrm{TT}$-algebras, chain Hamiltonians, potentials
for factorizable $S$-matrices, non-commutative spaces and link
invariants. Such aspects deserve further study, which is beyond
the scope of this paper. Some points are discussed in the
concluding remarks, including certain aspects of entanglements.

\section{Constructions (Even dimensions)}
\setcounter{equation}{0}

The $\left(2n\right)^2\times \left(2n\right)^2$ ($n\geq 1$)
$\widehat{R}$ matrices satisfying
\begin{equation}
\widehat{R}_{12}\widehat{R}_{23}\widehat{R}_{12}=\widehat{R}_{23}\widehat{R}_{12}\widehat{R}_{23}.
\end{equation}
are constructed in terms of the following operators which are
$\left(2n\right)\times \left(2n\right)$ matrices. Define
\begin{eqnarray}
&&I=\sum_{i=1}^n\left((ii)+(\overline{i}\overline{i})\right),
\qquad
J=\sum_{i=1}^n\left((-1)^{\bar{i}}(i\overline{i})+(-1)^{i}(\overline{i}i)\right),
\nonumber\\
&&K=\sum_{i=1}^n\left((i\overline{i})+(\overline{i}i)\right),
\qquad L=\sum_{i=1}^n\left((i\overline{i})-(\overline{i}i)\right),
\end{eqnarray}
where $\overline{i}=2n-i+1$ and $(ij)$ is the
$\left(2n\right)\times \left(2n\right)$ matrix with 1 for the
element (row $i$, column $j$) and zero elsewhere. The dimension
will be indicated (writing, say, $J_{(2n)}$ for $J$) when crucial
but not otherwise. For $n=1$, apart for the identity $I$
\begin{equation}
J=\left|\begin{array}{cc}
   0 & 1 \\
   -1 & 0\\
   \end{array}\right|,\qquad K=\left|\begin{array}{cc}
   0 & 1 \\
   1 & 0\\
   \end{array}\right|,\qquad L=J.
\end{equation}
The degeneracy is lifted ($J\neq L$) for $n>1$. For $n=2$,
\begin{equation}
J=\left|\begin{array}{cccc}
   0 & 0 & 0 & 1\\
   0 & 0 & -1 & 0\\
   0 & 1 & 0 & 0\\
   -1 & 0 & 0 & 0\\
   \end{array}\right|, \qquad K=\left|\begin{array}{cccc}
   0 & 0 & 0 & 1\\
   0 & 0 & 1 & 0\\
   0 & 1 & 0 & 0\\
   1 & 0 & 0 & 0\\
   \end{array}\right|,\qquad L=\left|\begin{array}{cccc}
   0 & 0 & 0 & 1\\
   0 & 0 & 1 & 0\\
   0 & -1 & 0 & 0\\
   -1 & 0 & 0 & 0\\
   \end{array}\right|.
 \end{equation}
Noting that $(-1)^{\overline{i}}=(-1)^{2n-i+1}=-(-1)^i$, one
obtains
\begin{eqnarray}
&&JK=-KJ=\sum_{i=1}^n\left((-1)^{\overline{i}}(ii)+(-1)^i(\overline{i}\overline{i})\right),
\nonumber\\
&&LK=-KL=\sum_{i=1}^n\left((ii)-(\overline{i}\overline{i})\right),\nonumber\\
&&JL=LJ=\sum_{i=1}^n\left((-1)^i(ii)-(-1)^{\overline{i}}(\overline{i}\overline{i})\right)
\end{eqnarray}
and
\begin{eqnarray}
&&-J^2=K^2=-L^2=I.
\end{eqnarray}
%(The $L_{ij}^{\pm}$ operators, to be defined later, will be quite
%distinct form $L$ without indices defined above. This $L$ will be
%used in diagonalization.)

The $\widehat{R}$ matrix, presented in sec. 1 (for $n=1$), is
generalized (for $n>1$) as follows:
\paragraph{Class I:}
\begin{equation}
\widehat{R}^{\pm 1}=\frac 1{\sqrt{2}}\left(I\otimes I\pm K\otimes
J\right)\end{equation}
satisfying unitarity
\begin{equation}
\widehat{R}^{+}\widehat{R}=I\otimes I=I_{(2n)^2}\end{equation} and
also
\begin{equation}
\widehat{R}^{2}=\sqrt{2}\widehat{R}-I\qquad \Rightarrow\qquad
\widehat{R}+\widehat{R}^{-1}=\sqrt{2}I.\end{equation} Using (2.5),
(2.6) one obtains
\begin{eqnarray}
&&\left(K\otimes J\otimes I\right)\left(I\otimes K\otimes
J\right)\left(K\otimes J\otimes I\right)=\left(I\otimes K\otimes
J\right),\nonumber\\
&& \left(I\otimes K\otimes J\right)\left(K\otimes J\otimes
I\right)\left(I\otimes K\otimes J\right)=\left(K\otimes J\otimes
I\right)\end{eqnarray} leading to
\begin{equation}
\widehat{R}_{12}\widehat{R}_{23}\widehat{R}_{12}=\frac
1{\sqrt{2}}\left(I\otimes K\otimes J+K\otimes J\otimes I\right)
=\widehat{R}_{23}\widehat{R}_{12}\widehat{R}_{23}.\end{equation}
Thus the braid equation is satisfied.
\paragraph{Baxterization:} Define $\left(z=\tanh\theta\right)$
\begin{eqnarray}
&&\widehat{R}\left(z\right)^{\pm 1}=\frac
1{\sqrt{1+z^2}}\left(I\otimes I\pm z K\otimes
J\right)\nonumber\\
&&\phantom{\widehat{R}\left(z\right)^{\pm 1}}=\frac
1{\sqrt{1+\left(\tanh\theta\right)^2}}\left(I\otimes I\pm
\tanh\theta K\otimes J\right)\equiv
\widehat{R}\left(\theta\right)^{\pm 1},\\
&&\widehat{R}\left(\pm 1\right)=\widehat{R}^{\pm 1}.\end{eqnarray}
They satisfy unitarity (for real $z$)\footnote{For real $\theta$
and $-\infty<\theta<+\infty$ one has the domain $-1<z<+1$. The
linearity in $z$ of $\sqrt{z^2+1}\widehat{R}\left(z\right)$ is
particularly helpful.}
\begin{equation}
\widehat{R}^{+}\left(z\right)\widehat{R}\left(z\right)=I\otimes I
\end{equation}
and
\begin{equation}
\widehat{R}\left(z\right)+\widehat{R}^{-1}\left(z\right)=\frac
2{\sqrt{1+z^2}}I\otimes I.
\end{equation}
Using (2.5), (2.6) and (2.10), we obtain
\begin{eqnarray}
&&\left(1+z^2\right)^{3/2}\widehat{R}_{12}\left(z\right)\widehat{R}_{23}\left(z''\right)\widehat{R}_{12}\left(z'\right)=
\left(1-zz'\right)\left(I\otimes I\otimes
I\right)+\nonumber\\
&&\phantom{\left(1+z^2\right)^{3/2}\widehat{R}_{12}\left(z\right)\widehat{R}_{23}\left(z''\right)
\widehat{R}_{12}\left(z'\right)=}z''\left(1+zz'\right)\left(I\otimes
K\otimes J\right)+\nonumber\\
&&\phantom{\left(1+z^2\right)^{3/2}\widehat{R}_{12}\left(z\right)\widehat{R}_{23}\left(z''\right)\widehat{R}_{12}
\left(z'\right)=}\left(z+z'\right)\left(K\otimes J\otimes
I\right)+\nonumber\\
&&\phantom{\left(1+z^2\right)^{3/2}\widehat{R}_{12}\left(z\right)\widehat{R}_{23}\left(z''\right)\widehat{R}_{12}
\left(z'\right)=}z''\left(z'-z\right)\left(K\otimes KJ\otimes I\right),\\
&&\left(1+z^2\right)^{3/2}\widehat{R}_{23}\left(z'\right)\widehat{R}_{12}\left(z''\right)\widehat{R}_{23}\left(z\right)=
\left(1-zz'\right)\left(I\otimes I\otimes
I\right)+\nonumber\\
&&\phantom{\left(1+z^2\right)^{3/2}\widehat{R}_{12}\left(z\right)\widehat{R}_{23}\left(z''\right)\widehat{R}_{12}
\left(z'\right)=}\left(z+z'\right)\left(I\otimes K\otimes J\right)+ \nonumber\\
&&\phantom{\left(1+z^2\right)^{3/2}\widehat{R}_{12}\left(z\right)\widehat{R}_{23}\left(z''\right)\widehat{R}_{12}
\left(z'\right)=}z''\left(1+zz'\right)\left(K\otimes J\otimes
I\right)+\nonumber\\
&&\phantom{\left(1+z^2\right)^{3/2}\widehat{R}_{12}\left(z\right)\widehat{R}_{23}\left(z''\right)\widehat{R}_{12}
\left(z'\right)=}z''\left(z'-z\right)\left(K\otimes KJ\otimes
I\right).\end{eqnarray} Setting
\begin{equation}
z''=\frac{z+z'}{1+zz'}=\frac{\tanh\theta+\tanh\theta'}{1+\tanh\theta\tanh\theta'}=\tanh\left(\theta+\theta'\right).
\end{equation}
one obtains the baxterized braid equation (spectral parameter
dependent generalization of (2.11))
\begin{equation}
\widehat{R}_{12}\left(z\right)\widehat{R}_{23}\left(z''\right)\widehat{R}_{12}\left(z'\right)=
\widehat{R}_{23}\left(z'\right)\widehat{R}_{12}\left(z''\right)\widehat{R}_{23}\left(z\right)
\end{equation}
or
\begin{equation}
\widehat{R}_{12}\left(\theta\right)\widehat{R}_{23}\left(\theta+\theta'\right)\widehat{R}_{12}\left(\theta'\right)=
\widehat{R}_{23}\left(\theta'\right)\widehat{R}_{12}\left(\theta+\theta'\right)\widehat{R}_{23}\left(\theta\right).
\end{equation}

\paragraph{Class II:} We briefly state that $P\widehat{R}P$ of
(1.1) is generalized as
\begin{eqnarray}
&&\widehat{R}^{\pm 1}=\frac 1{\sqrt{2}}\left(I\otimes I\pm
J\otimes K\right),\nonumber\\
&&\widehat{R}\left(z\right)^{\pm 1}=\frac
1{\sqrt{1+z^2}}\left(I\otimes I\pm z J\otimes K\right).
\end{eqnarray}
The braid equation and its baxterization are verified following
closely the steps for the preceding case. Unitarity is preserved.
Defining the $\left(2n\right)^2\times\left(2n\right)^2$
permutation matrix
\begin{equation}
P=\sum_{a,b=1}^{2n}(ab)\otimes (ba),
\end{equation} we obtain
\begin{equation}
\widehat{R}_{(II)}\left(z\right)=P\widehat{R}_{(I)}\left(z\right)P,
\end{equation}
where $\widehat{R}_{(I)}\left(z\right)$,
$\widehat{R}_{(II)}\left(z\right)$ correspond to
$\widehat{R}\left(z\right)$ for class (I) and class (II)
respectively.

\section{Projectors}
\setcounter{equation}{0}

We start with the spectral resolution for class (I). Define
\begin{equation}
P_{\pm}=\frac 12 \left(I\otimes I\pm \mathrm{i}K\otimes J\right).
\end{equation}
Using (2.5) and (2.6) one obtains $\left(\epsilon,\epsilon'=\pm
1\right)$
\begin{eqnarray}
&&P_{\epsilon}P_{\epsilon'}=\frac 14 \left(I\otimes
I+\mathrm{i}\left(\epsilon+\epsilon'\right)K\otimes
J+\epsilon\epsilon' I\otimes I\right)\nonumber\\
&&\phantom{P_{\epsilon}P_{\epsilon'}}=\frac 14
\left(\left(1+\varepsilon\varepsilon'\right)I\otimes
I+\mathrm{i}\left(\epsilon+\epsilon'\right)K\otimes
J\right)=P_{\epsilon}\delta_{\epsilon,\epsilon'}.
\end{eqnarray}
On such a basis
\begin{eqnarray}
&&\widehat{R}^{\pm 1}=\frac 1{\sqrt{2}} \left(\left(1\mp
\mathrm{i}\right)P_++\left(1\pm \mathrm{i}\right)P_-\right), \\
&&\widehat{R}\left(z\right)^{\pm 1}=\frac 1{\sqrt{1+z^2}}
\left(\left(1\mp \mathrm{i}z\right) P_++\left(1\pm
\mathrm{i}z\right)P_-\right).
\end{eqnarray}
Such a basis with {\it complex} projectors for {\it real} $R$ and
$R\left(z\right)$ was already implemented in our previous study of
$\mathrm{SO3}$ for $n=1$ (see refs. \cite{R1,R2}). For $n=2$,
(3.1) implies that
\begin{eqnarray}
&&P_{\pm 1}=\frac 12\left|\begin{array}{cccc}
   I & 0 & 0 & \pm \mathrm{i}J\\
   0 & I & \pm \mathrm{i}J & 0\\
   0 & \pm \mathrm{i}J & I & 0\\
   \pm \mathrm{i}J & 0 & 0 & I\\
   \end{array}\right|=\pm \frac {\mathrm{i}}{\sqrt{2}}\left(\widehat{R}-
   \frac{\left(1\pm \mathrm{i}\right)}{\sqrt{2}}I\otimes I\right),\\
&&P_++P_-=I\otimes I.
\end{eqnarray}
The projectors play basic roles in the construction of
non-commutative spaces associated to $\widehat{R}$. A parallel
treatment of projectors for class (II) can evidently be carried
through. It will not be presented explicitly.

\section{Diagonalization, block-diagonalization and a non-equivalence}
\setcounter{equation}{0}

Here the term {\it non-equivalence} refers to non-conservation of
the braid equation. Define
\begin{equation}
M^{\pm 1}=\frac 1{\sqrt{2}}\left(I\otimes I\pm \mathrm{i}L\otimes
J\right)
\end{equation}
giving, say, for $n=2$
\begin{equation}
M^{\pm 1}=\frac 1{\sqrt{2}}\left|\begin{array}{cccc}
   I & 0 & 0 & \pm \mathrm{i}J\\
   0 & I & \pm \mathrm{i}J & 0\\
   0 & \mp \mathrm{i}J & I & 0\\
   \mp \mathrm{i}J & 0 & 0 & I\\
   \end{array}\right|,
\end{equation}
one obtains (for $\epsilon=\pm$)
\begin{equation}
MP_\epsilon M^{-1}=\frac 12\left(I+\epsilon LK\right)\otimes I.
\end{equation}
For $n=2$, for example,
\begin{eqnarray}
&&MP_+M^{-1}=\left|\begin{array}{cccc}
   I & 0 & 0 & 0\\
   0 & I & 0 & 0\\
   0 & 0 & 0 & 0\\
   0 & 0 & 0 & 0\\
   \end{array}\right|,
   \qquad MP_-M^{-1}=\left|\begin{array}{cccc}
   0 & 0 & 0 & 0\\
   0 & 0 & 0 & 0\\
   0 & 0 & I & 0\\
   0 & 0 & 0 & I\\
   \end{array}\right|,\\
   && M\widehat{R}\left(z\right)M^{-1}=\frac 1{\sqrt{1+z^2}}\left|\begin{array}{cccc}
   (1-\mathrm{i}z)I & 0 & 0 & 0\\
   0 & (1-\mathrm{i}z)I & 0 & 0\\
   0 & 0 & (1+\mathrm{i}z)I & 0\\
   0 & 0 & 0 & (1+\mathrm{i}z)I\\
   \end{array}\right|.
\end{eqnarray}
For $z=\pm 1$ one obtains the results for $\widehat{R}^{\pm 1}$
respectively. In (4.5), for $n=2$, $I\equiv
I_{(4)}=\left|\begin{array}{cc}
   I_{(2)} & 0 \\
   0 & I_{(2)} \\
   \end{array}\right|$,
$I_{2n}$ being the identity matrix of $(2n)\times (2n)$
dimensions. Similarly in (4.2), $J\equiv
J_{(4)}=\left|\begin{array}{cc}
   0 & J_{(2)} \\
   J_{(2)} & 0 \\
   \end{array}\right|$.
Starting with, for $2n=2$, and setting $z=1$ for simplicity,
\begin{equation}
\widehat{R}=\frac 1{\sqrt{2}}\left|\begin{array}{cc}
   I_{(2)} & J_{(2)} \\
   J_{(2)} & I_{(2)} \\
   \end{array}\right|\equiv \widehat{R}_{(2)}
\end{equation}
can be diagonalized by conjugating with
\begin{equation}
M_{(2)}^{\pm 1}=\frac 1{\sqrt{2}}\left|\begin{array}{cc}
   I_{(2)} & \pm \mathrm{i}J_{(2)} \\
   \mp \mathrm{i}J_{(2)} & I_{(2)} \\
   \end{array}\right|
\end{equation}
giving
\begin{equation}
M_{(2)}R_{(2)}M_{(2)}^{\pm 1}=\frac
1{\sqrt{1+z^2}}\left|\begin{array}{cc}
   (1-\mathrm{i}z)I_2 & 0 \\
   0 & (1+\mathrm{i}z)I_2 \\
   \end{array}\right|.
\end{equation}
Such diagonalizations indicates clearly how $R_{(4)}$, i.e.
$\widehat{R}$ (for $n=2$) can be block diagonalized into a direct
sum of 4 successive $R_{(2)}$, namely (suppressing the argument
$z$ for simplicity)
\begin{equation}
VR_{(4)}V^{-1}=\left|\begin{array}{cccc}
   R_{(2)} & 0 & 0 & 0\\
   0 & R_{(2)} & 0 & 0 \\
   0 & 0 & R_{(2)} & 0\\
   0 & 0 & 0 & R_{(2)}\\
   \end{array}\right|\equiv \widehat{R}'.
\end{equation}
First one permutes the $2\times 2$ blocks of (4.5) by conjugating
with
\begin{equation}
U=U^{-1}=\left|\begin{array}{cccccccc}
   I_{(2)} & 0 & 0 & 0 &0 & 0 &0 & 0 \\
   0 & 0 & 0 & 0 &0 & 0 &I_{(2)} & 0 \\
   0 & 0 & I_{(2)} & 0 &0 & 0 &0 & 0 \\
   0 & 0 & 0 & 0 &I_{(2)} & 0 &0 & 0 \\
   0 & 0 & 0 & I_{(2)} &0 & 0 &0 & 0 \\
   0 & 0 & 0 & 0 &0 & I_{(2)} &0 & 0 \\
   0 & I_{(2)} & 0 & 0 &0 & 0 &0 & 0 \\
   0 & 0 & 0 & 0 &0 & 0 &0 & I_{(2)} \\
   \end{array}\right|.
\end{equation}
Then conjugate this back to (4.9) by blocks
$\left(M^{-1}_{(2)},M_{(2)}\right)$ in block diagonal form (4
blocks of $M^{\mp 1}_{(2)}$). The combined conjugation gives $V$
and (4.9). But one can easily see that $V$ is not of the form of a
tensor product of some matrix $Y$ (of $4\times 4$ dimension) i.e.
\begin{equation}
V\neq Y\otimes Y,
\end{equation}
where conjugation by invertible $Y\otimes Y$ conserves braid
property. This has for consequence that (4.9) does not satisfy the
braid equation. Direct computation gives (for $n=2$, for example)
\begin{equation}
\widehat{R}_{12}'\widehat{R}_{23}'\widehat{R}_{12}'-\widehat{R}_{23}\widehat{R}_{12}'\widehat{R}_{23}'=
\left(I_4\otimes \widehat{R}_{(2)}\otimes I_4-I_4\otimes
I_4\otimes \widehat{R}_{(2)}\right)\neq 0.
\end{equation}
(One uses $\widehat{R}_{12}'=I_4\otimes \widehat{R}_{(2)}\otimes
I_4$, $\widehat{R}_{23}'=I_4\otimes I_4\otimes \widehat{R}_{(2)}$;
$\widehat{R}_{(2)}^2=\sqrt{2}\widehat{R}_{(2)}-I$.) Such
considerations can easily be generalized to $n>2$. They show quite
explicitly that our generalizations (2.7) (and similarly (2.21))
are intrinsically non-equivalent to direct sums of the $n=1$ (i.e.
$4\times 4$) blocks, which do not conserve the braid property.
This holds though the forms can be related via a conjugation (by
$V$).

\paragraph{Further possibilities interchanging roles of $J$ and
$L$:} From the structure of the algebra (2.5), (2.6), it is
evident hat one can replace $\left(I,J,K,L\right)$. in the
preceding developments by $\left(I,L,K,J\right)$ respectively
retaining the essential results for $n>1$. For $n=1$ the two sets
coincide (since $J=L$). The two treatments can be related through
suitable permutations of rows and columns. We will not present
this aspect explicitly. But since $L$ appears in the diagonalizer
$M$ of (4.1), the full scope of the operator $L$ is worth noting.
Starting with
\begin{equation}
\widehat{R}\left(z\right)^{\pm 1}=\frac
1{\sqrt{1+z^2}}\left(I\otimes I\pm zK\otimes L\right)
\end{equation}
or
\begin{equation}
\widehat{R}\left(z\right)^{\pm 1}=\frac
1{\sqrt{1+z^2}}\left(I\otimes I\pm zL\otimes K\right)
\end{equation}
one again obtains, analogously to (2.16) the braid equation (with
$L$ replacing $J$). For $n=1$, (4.13) and (4.14) coincides with
(2.12) and (2.21) respectively. The different formulations are
related through permutations of appropriate rows an columns. One
retains a symmetric diagonal (unity) and a antisymmetric
anti-diagonal. Such strong constraints conserve unitarity and
braid property.

\section{Odd dimensions (A class of complex unitary braid matrices)}
\setcounter{equation}{0}

In previous papers \cite{R6,R8} $(2n+1)^2\times (2n+1)^2$
dimensional braid matrices have been constructed and studied for
$n>1$. They were obtained by implementing a {\it nested sequence}
of projectors defined already before (see Ref. \cite{R6}). In
these sources they were studied as real, symmetric braid matrices
with multiple parameters and already Baxterized \cite{R6}. For all
parameters real (as also the spectral parameter $\theta$) they
satisfy
\begin{equation}
\widehat{R}^+\left(\theta\right)=\widehat{R}\left(\theta\right),\qquad
\widehat{R}\left(-\theta\right)=\widehat{R}\left(\theta\right)^{-1}.
\end{equation}
Here we note that:
\begin{enumerate}
    \item for all parameters pure imaginary and $\theta$ real, or
    alternatively
    \item for all parameters real and $\theta$ pure imaginary
\end{enumerate}
they become unitary, i.e.
\begin{equation}
\widehat{R}^+\left(\theta\right)=\widehat{R}\left(-\theta\right)=\widehat{R}\left(\theta\right)^{-1}.
\end{equation}
This happens due to the special structure of this class. It is
sufficient to illustrate this for $9\times 9$ matrix ($n=1$). This
involves six parameters \cite{R6} $\left(m_{11}^{\pm},
m_{12}^{\pm},m_{21}^{\pm}\right)$. Making them explicitly pure
imaginary as $m_{ij}^{\pm}\rightarrow \mathrm{i}m_{ij}^{\pm}$ with
real $m$'s on the right and defining
\begin{eqnarray}
&&a_{\pm}=\frac
12\left(\exp\left(\mathrm{i}m_{11}^{(+)}\theta\right)\pm
\exp\left(\mathrm{i}m_{11}^{(-)}\theta\right)\right), \nonumber\\
&& b_{\pm}=\frac
12\left(\exp\left(\mathrm{i}m_{12}^{(+)}\theta\right)\pm
\exp\left(\mathrm{i}m_{12}^{(-)}\theta\right)\right),\nonumber\\
&& c_{\pm}=\frac
12\left(\exp\left(\mathrm{i}m_{21}^{(+)}\theta\right)\pm
\exp\left(\mathrm{i}m_{21}^{(-)}\theta\right)\right)
\end{eqnarray}
one obtains
\begin{equation}
a_{+}\overline{a_+}+a_{-}\overline{a_-}=1,\qquad
a_{+}\overline{a_-}+a_{-}\overline{a_+}=0
\end{equation}
and so on. Now it is easy to see that
\begin{equation}
\widehat{R}\left(\theta\right)=\left|\begin{array}{ccccccccc}
   a_+ & 0 & 0 & 0 & 0 & 0 & 0 & 0 & a_-\\
   0 & b_+ & 0 & 0 & 0 & 0 & 0 & b_- & 0\\
   0 & 0 & a_+ & 0 & 0 & 0 & a_- & 0 & 0\\
   0 & 0 & 0 & c_+ & 0 & c_- & 0 & 0 & 0\\
   0 & 0 & 0 & 0 & 1 & 0 & 0 & 0 & 0\\
   0 & 0 & 0 & c_- & 0 & c_+ & 0 & 0 & 0\\
   0 & 0 & a_- & 0 & 0 & 0 & a+ & 0 & 0\\
   0 & b_- & 0 & 0 & 0 & 0 & 0 & b_+ & 0\\
   a_- & 0 & 0 & 0 & 0 & 0 & 0 & 0 & a_+\\
   \end{array}\right|
\end{equation}
satisfies
\begin{equation}
\widehat{R}\left(\theta\right)^+\widehat{R}\left(\theta\right)=\widehat{R}\left(-\theta\right)
\widehat{R}\left(\theta\right)=I.
\end{equation}
Evidently, the braid equation is still satisfied since that does
not depend on the reality condition. The generalization for $n>1$
is trivial. Thus we obtain a class of complex, unitary braid
matrices for odd dimensions. Such a class is, however, well
defined only in the Baxterized form. The limits of {\it infinite
rapidity} ($\theta\rightarrow\pm \infty$ or $z\rightarrow \pm 1$
in the preceding even dim. ones) give here oscillating
exponentials. Since
\begin{equation}
\exp\left(\mathrm{i}m_{il}^{\pm}\theta\right)=\exp\left[{\mathrm{i}
m_{il}^{\pm}\left(\theta+\frac{2\pi}{m_{il}^{\pm}}\right)}\right]
\end{equation}
if all the m's are commensurate (with rational ratios) there will
be an overall common period for all the parameters  as $\theta$
varies. Then $\widehat{R}\left(\theta\right)$ is periodic in
$\theta$. But if at least two $m$'s are incommensurate
$\widehat{R}\left(\theta\right)$ is quasi-periodic. Such aspects
might be worth study.

Further exploration of complex unitary braid matrices is beyond
the scope of this paper. In the following sections only real, even
dimensional unitary braid matrices are studied. See, however, the
remarks in conclusion.

\section{$\widehat{R}\mathrm{LL}$- and $\widehat{R}\mathrm{TT}$-algebras}
\setcounter{equation}{0}

For classes I and II the $\mathrm{LL}$- and $\mathrm{TT}$-algebras
are simply interchanged due to the relation
\begin{equation}
P\widehat{R}_{(I)}\left(z\right)P=\frac 1{\sqrt{1+z^2}}
P\left(I\otimes I +zK\otimes J\right)P=\frac 1{\sqrt{1+z^2}}
\left(I\otimes I +zJ\otimes K\right)= \widehat{R}_{(II)}
\left(z\right)
\end{equation}
and the fact that the fundamental blocks of $\mathrm{L}$ and
$\mathrm{T}$ are obtained from
\begin{equation}
\mathrm{L}\left(z\right)=\widehat{R}\left(z\right)P,\qquad
\mathrm{T}\left(z\right)=P\widehat{R}\left(z\right)
\end{equation}
leading to
\begin{eqnarray}
&&\mathrm{L}_{(II)}\left(z\right)=\widehat{R}_{(II)}\left(z\right)P=
P\widehat{R}_{(I)}\left(z\right)=\mathrm{T}_{(I)}\left(z\right),\\
&&\mathrm{T}_{(II)}\left(z\right)=P\widehat{R}_{(II)}\left(z\right)=
\widehat{R}_{(I)}\left(z\right)P=\mathrm{L}_{(I)}\left(z\right).
\end{eqnarray}

To obtain higher order representations one implements the same
coproduct prescriptions for $\mathrm{L}$ and $\mathrm{T}$. Hence
the correspondence is maintained. We study below only the case I,
suppressing the index ($\widehat{R}_{(I)}\rightarrow \widehat{R}$
and so on). Corresponding to
\begin{equation}
\widehat{R}\left(z\right)=\frac 1{\sqrt{2\left(1+z^2\right)}}
\left(\left(1+z\right)\widehat{R}+(1-z)\widehat{R}^{-1}\right)
\end{equation}
one can define (simplifying the external factor irrelevant for our
purposes)
\begin{eqnarray}
&&\mathrm{L}\left(z\right)=\frac 12
\left(\left(1+z\right)\mathrm{L}^++(1-z)\mathrm{L}^{-}\right)=
\frac {e^\theta \mathrm{L}^++e^{-\theta}\mathrm{L}^{-}}{e^\theta+e^{-\theta}} ,\\
&&\mathrm{L}\left(\pm 1\right)=\mathrm{L}^{\pm}.
\end{eqnarray}
$(z=\tanh\theta$). The single constraint (with
$\mathrm{L}_1=\mathrm{L}\otimes I$, $\mathrm{L}_2=I\otimes
\mathrm{L}$)
\begin{equation}
\widehat{R}\left(\theta-\theta'\right)\mathrm{L}_2\left(\theta\right)\mathrm{L}_1\left(\theta'\right)=
\mathrm{L}_2\left(\theta'\right)\mathrm{L}_1\left(\theta\right)\widehat{R}\left(\theta-\theta'\right)
\end{equation}
can be shown to imply all the three FRT relations \cite{R2}
\begin{equation}
\widehat{R}\mathrm{L}_2^{\epsilon}\mathrm{L}_1^{\epsilon'}=
\mathrm{L}_2^{\epsilon'}\mathrm{L}_1^{\epsilon}\widehat{R},
\end{equation}
where
$\left(\epsilon,\epsilon'\right)=\left(+,+\right),\left(-,-\right),\left(+,-\right)$
respectively.

Constructions such as (6.6) are only possible when
$\widehat{R}\left(z\right)$ satisfies a quadratic constraint (i.e.
$\widehat{R}\left(z\right)$, $\widehat{R}\left(z\right)^{-1}$ a
linear one). In all our constructions (2.15) and (6.5) are
guaranteed and hence also (6.8). A quadratic constraint and hence
(6.6) can be shown to permit two distinct type of coproducts which
coincide for $z=\pm 1$ but are inequivalent for ($-1<z<1$)
\cite{R2}. In this paper we will, for brevity, restrict our study
to the standard prescription,
\begin{equation}
\mathrm{L}_{ij}^{(r+1)}\left(z\right)=\sum_k\mathrm{L}_{ik}^{(1)}\left(z\right)\otimes
\mathrm{L}_{kj}^{(r)}\left(z\right),
\end{equation}
where $\mathrm{L}_{ik}^{(1)}\left(z\right)$ are obtained from
(6.2) as follows:
\begin{eqnarray}
&&\mathrm{L}^{(1)}\left(z\right)\equiv\mathrm{L}\left(z\right)=\widehat{R}\left(z\right)P=\frac
1{\sqrt{z^2+1}}\sum_{i,j=1}^{2n}\left((ij)\otimes
(ji)+z(-1)^{\overline{j}}(i\overline{j})\otimes
(j\overline{i})\right)\nonumber\\
&&\phantom{\mathrm{L}^{(1)}\left(z\right)}\equiv
\left|\begin{array}{cccc}
\mathrm{L}_{11}^{(1)}\left(z\right)&\mathrm{L}_{12}^{(1)}\left(z\right)& \cdots & \mathrm{L}_{1,2n}^{(1)}\left(z\right)\\
\mathrm{L}_{21}^{(1)}\left(z\right)&\mathrm{L}_{22}^{(1)}\left(z\right)& \cdots & \mathrm{L}_{2,2n}^{(1)}\left(z\right)\\
\vdots&\vdots& \ddots & \vdots\\
\mathrm{L}_{2n,1}^{(1)}\left(z\right)&\mathrm{L}_{2n,2}^{(1)}\left(z\right)&
\cdots &
\mathrm{L}_{2n,2n}^{(1)}\left(z\right)\\\end{array}\right|.
\end{eqnarray}
For $n=2$, for example, (6.10) gives (suppressing for simplicity
the argument $z$ for each $\mathrm{L}_{ij}$ and dropping the
overall factor)
\begin{eqnarray}
&&\mathrm{L}_{1j}^{(r+1)}=\left|\begin{array}{cccc}
   \mathrm{L}_{1j}^{(r)} & 0 & 0 & z\mathrm{L}_{4j}^{(r)} \\
   \mathrm{L}_{2j}^{(r)} & 0 & 0 & -z\mathrm{L}_{3j}^{(r)} \\
   \mathrm{L}_{3j}^{(r)} & 0 & 0 & z\mathrm{L}_{2j}^{(r)} \\
   \mathrm{L}_{4j}^{(r)} & 0 & 0 & -z\mathrm{L}_{1j}^{(r)} \\
   \end{array}\right|,\qquad
   \mathrm{L}_{2j}^{(r+1)}=\left|\begin{array}{cccc}
   0 &\mathrm{L}_{1j}^{(r)} & z\mathrm{L}_{4j}^{(r)} & 0\\
   0 & \mathrm{L}_{2j}^{(r)} & -z\mathrm{L}_{3j}^{(r)} & 0\\
   0 & \mathrm{L}_{3j}^{(r)} & z\mathrm{L}_{2j}^{(r)} & 0 \\
   0 & \mathrm{L}_{4j}^{(r)} & -z\mathrm{L}_{1j}^{(r)} & 0\\
   \end{array}\right|,\nonumber\\
&&\mathrm{L}_{3j}^{(r+1)}=\left|\begin{array}{cccc}
   0 & z\mathrm{L}_{4j}^{(r)} &   \mathrm{L}_{1j}^{(r)} & 0\\
   0 & -z\mathrm{L}_{3j}^{(r)}   &  \mathrm{L}_{2j}^{(r)} & 0\\
   0 & z\mathrm{L}_{2j}^{(r)}  &   \mathrm{L}_{3j}^{(r)} & 0\\
   0 & -z\mathrm{L}_{1j}^{(r)}  & \mathrm{L}_{4j}^{(r)} & 0\\
   \end{array}\right|,\qquad
\mathrm{L}_{4j}^{(r+1)}=\left|\begin{array}{cccc}
   z\mathrm{L}_{4j}^{(r)} & 0 & 0 &  \mathrm{L}_{1j}^{(r)}\\
   -z\mathrm{L}_{3j}^{(r)} & 0 & 0 &  \mathrm{L}_{2j}^{(r)}\\
  z\mathrm{L}_{2j}^{(r)}  & 0 & 0 &  \mathrm{L}_{3j}^{(r)}\\
  -z\mathrm{L}_{1j}^{(r)}  & 0 & 0 &  \mathrm{L}_{4j}^{(r)}\\
   \end{array}\right|,
\end{eqnarray}
where $j=1,2,3,4$. Setting now
$\mathrm{L}_{ij}^{(0)}=\delta_{ij}$, we obtain
\begin{equation}
\mathrm{L}_{11}^{(1)}=\left|\begin{array}{cccc}
   1 & 0 & 0 & 0 \\
   0 & 0 & 0 & 0 \\
   0 & 0 & 0 & 0 \\
   0 & 0 & 0 & -z \\
   \end{array}\right|,\qquad
   \mathrm{L}_{12}^{(1)}=\left|\begin{array}{cccc}
   0 & 0 & 0 & 0 \\
   1 & 0 & 0 & 0 \\
   0 & 0 & 0 & z \\
   0 & 0 & 0 & 0 \\
   \end{array}\right|,\qquad \hbox{etc.}
\end{equation}

Define
\begin{equation}
\mathrm{L}^{(r+1)}=\sum_{i=1}^{2n}\mathrm{L}_{ii}^{(r+1)}.
\end{equation}
When $n=2$, for example, we have
\begin{equation}
\hbox{\verb"Tr"}\left(\mathrm{L}^{(r+1)}\right)=\sum_{i=1}^{4}\hbox{\verb"Tr"}\left(\mathrm{L}_{ii}^{(r+1)}\right)
=\hbox{\verb"Tr"}\left(\left(1-z\right)\left(\mathrm{L}_{11}^{(r)}+\mathrm{L}_{33}^{(r)}\right)+
\left(1+z\right)\left(\mathrm{L}_{22}^{(r)}+\mathrm{L}_{44}^{(r)}\right)\right).
\end{equation}
Observing that
\begin{eqnarray}
&&\hbox{\verb"Tr"}\left(\mathrm{L}_{11}^{(r)}\right)=\left(1-z\right)\hbox{\verb"Tr"}\left(\mathrm{L}_{11}^{(r-1)}\right),\qquad
\hbox{\verb"Tr"}\left(\mathrm{L}_{22}^{(r)}\right)=\left(1+z\right)\hbox{\verb"Tr"}\left(\mathrm{L}_{22}^{(r-1)}\right),\nonumber\\
&&\hbox{\verb"Tr"}\left(\mathrm{L}_{33}^{(r)}\right)=\left(1-z\right)\hbox{\verb"Tr"}\left(\mathrm{L}_{33}^{(r-1)}\right),\qquad
\hbox{\verb"Tr"}\left(\mathrm{L}_{44}^{(r)}\right)=\left(1+z\right)\hbox{\verb"Tr"}\left(\mathrm{L}_{44}^{(r-1)}\right),
\end{eqnarray}
we deduce
\begin{equation}
\hbox{\verb"Tr"}\left(\mathrm{L}_{ii}^{(r)}\right)=\left(1+(-1)^iz\right)^r,
\qquad i=1,2,3,4.
\end{equation}
Using (6.17), we finally obtain
\begin{equation}
\hbox{\verb"Tr"}\left(\mathrm{L}^{(r)}\right)=2\left(\left(1+z\right)^r+\left(1-z\right)^r\right).
\end{equation}

The $\widehat{R}\mathrm{TT}$ constraints are conveniently written
as
\begin{equation}
\widehat{R}\left(\theta-\theta'\right)\mathrm{T}\left(\theta\right)\otimes
\mathrm{T}\left(\theta'\right)=\mathrm{T}\left(\theta'\right)\otimes
\mathrm{T}\left(\theta\right)\widehat{R}\left(\theta-\theta'\right).
\end{equation}
Starting from (6.2) and a coproduct prescription parallel to
(6.10) one obtains analogously, for example, when $n=2$
(suppressing again arguments $z$)
\begin{eqnarray}
&&\mathrm{T}_{1j}^{(r+1)}=\left|\begin{array}{cccc}
   \mathrm{T}_{1j}^{(r)} & 0 & 0 & z\mathrm{T}_{4j}^{(r)} \\
   \mathrm{T}_{2j}^{(r)} & 0 & 0 & z\mathrm{T}_{3j}^{(r)} \\
   \mathrm{T}_{3j}^{(r)} & 0 & 0 & z\mathrm{T}_{2j}^{(r)} \\
   \mathrm{T}_{4j}^{(r)} & 0 & 0 & z\mathrm{T}_{1j}^{(r)} \\
   \end{array}\right|,\qquad
   \mathrm{T}_{2j}^{(r+1)}=\left|\begin{array}{cccc}
   0 &\mathrm{T}_{1j}^{(r)} & -z\mathrm{T}_{4j}^{(r)} & 0\\
   0 & \mathrm{T}_{2j}^{(r)} & -z\mathrm{T}_{3j}^{(r)} & 0\\
   0 & \mathrm{T}_{3j}^{(r)} & -z\mathrm{T}_{2j}^{(r)} & 0 \\
   0 & \mathrm{T}_{4j}^{(r)} & -z\mathrm{T}_{1j}^{(r)} & 0\\
   \end{array}\right|,\nonumber\\
&&\mathrm{T}_{3j}^{(r+1)}=\left|\begin{array}{cccc}
   0 & z\mathrm{T}_{4j}^{(r)} &  \mathrm{T}_{1j}^{(r)} & 0\\
   0 & z\mathrm{T}_{3j}^{(r)} &  \mathrm{T}_{2j}^{(r)} & 0\\
   0 & z\mathrm{T}_{2j}^{(r)} &  \mathrm{T}_{3j}^{(r)} & 0\\
   0 & z\mathrm{T}_{1j}^{(r)} &  \mathrm{T}_{4j}^{(r)} & 0\\
   \end{array}\right|,\qquad
\mathrm{T}_{4j}^{(r+1)}=\left|\begin{array}{cccc}
   -z\mathrm{T}_{4j}^{(r)} & 0 & 0 &  \mathrm{T}_{1j}^{(r)}\\
   -z\mathrm{T}_{3j}^{(r)} & 0 & 0 &  \mathrm{T}_{2j}^{(r)}\\
   -z\mathrm{T}_{2j}^{(r)} & 0 & 0 &  \mathrm{T}_{3j}^{(r)}\\
   -z\mathrm{T}_{1j}^{(r)} & 0 & 0 &  \mathrm{T}_{4j}^{(r)}\\
   \end{array}\right|.
\end{eqnarray}
By setting again $\mathrm{T}_{ij}^{(0)}=\delta_{ij}$, we obtain
\begin{equation}
\mathrm{T}_{11}=\left|\begin{array}{cccc}
   1 & 0 & 0 & 0 \\
   0 & 0 & 0 & 0 \\
   0 & 0 & 0 & 0 \\
   0 & 0 & 0 & z \\
   \end{array}\right|,\qquad
   \mathrm{T}_{12}=\left|\begin{array}{cccc}
   0 & 0 & 0 & 0 \\
   1 & 0 & 0 & 0 \\
   0 & 0 & 0 & z \\
   0 & 0 & 0 & 0 \\
   \end{array}\right|, \qquad \hbox{etc.}
\end{equation}
Similarly, we obtain, for $n=2$
\begin{eqnarray}
&&\hbox{Tr}\left(\mathrm{T}^{(r+1)}\right)=\sum_{i=1}^{4}\hbox{\verb"Tr"}\left(\mathrm{T}_{ii}^{(r+1)}\right)
\nonumber\\
&&\phantom{\hbox{Tr}\left(\mathrm{T}^{(r+1)}\right)}=\hbox{\verb"Tr"}\left(\left(1+z\right)\left(\mathrm{T}_{11}^{(r)}
+\mathrm{T}_{33}^{(r)}\right)+\left(1-z\right)\left(\mathrm{T}_{22}^{(r)}+\mathrm{T}_{44}^{(r)}\right)\right),\nonumber\\
&&\phantom{\hbox{Tr}\left(\mathrm{T}^{(r+1)}\right)}=\hbox{\verb"Tr"}\left(\left(1+z\right)^2\left(\mathrm{T}_{11}^{(r-1)}+
\mathrm{T}_{33}^{(r-1)}\right)+
\left(1-z\right)^2\left(\mathrm{t}_{22}^{(r-1)}+\mathrm{T}_{44}^{(r-1)}\right)\right),\nonumber\\
&&\phantom{\hbox{Tr}\left(\mathrm{T}^{(r+1)}\right)=}\vdots\nonumber\\
&&\hbox{Tr}\left(\mathrm{T}^{(r+1)}\right)=2\left(\left(1+z\right)^r+\left(1-z\right)^r\right).
\end{eqnarray}

The $\mathrm{L}$-algebra, for $n=1$, has been studied extensively
in a previous paper \cite{R2} which provides further references.
Here we have presented briefly certain basic features of
$\mathrm{L}_{ij}^{\pm}$ and $\mathrm{T}_{ij}$ for all orders
$(r)$. Their further study concerning representations will not be
undertaken in this paper. One application of the
$\mathrm{T}$-matrices concerns the Hamiltonians encoding the
evolution of the states in base space, whether the states are
entangled or not. Hamiltonians are briefly presented in the next
section.

Denoting the diagonal matrix (sec. 4)
($M\widehat{R}\left(z''\right)M^{-1}$) as $D\left(z''\right)$ for
all $n$ and $z''=\tanh\left(z-z'\right)$ one reduces (6.19) to
\begin{equation}
D\left(z''\right)\left(M\mathrm{T}\left(z\right)\otimes
\mathrm{T}\left(z'\right)M^{-1}\right)=\left(M\mathrm{T}\left(z'\right)\otimes
\mathrm{T}\left(z\right)M^{-1}\right)D\left(z''\right).
\end{equation}
This expresses the $\mathrm{TT}'$-algebra in terms of sums of
terms such that permutation of $z$ and $z'$ has extremely simple
consequences. The $\mathrm{LL}'$-algebra can be treated similarly.

The $\mathrm{TT}$- and $\mathrm{LL}$-algebras and their
representations will be studied in detail in our follow-up
paper(s) \cite{ACDM} along the lines of \cite{R1,R2,ACDM1},
namely, presentation of the $\mathrm{TT}$-algebras as matrix
bialgebras, construction of their dual bialgebras, presentation,
of the $\mathrm{LL}$-algebras as FRT-duals of the
$\mathrm{TT}$-algebras, development of the representation theory
of all mentioned algebras.

\section{Hamiltonians}
\setcounter{equation}{0}

Sections on Hamiltonians in previous papers \cite{R8,R9} cite
various references. Here we briefly present some generic features
arising from the structure of our unitary braid matrices. The
Hamiltonian, of order $r$, is defined in terms of the transfer
matrix $\mathrm{T}^{(r)}\left(\theta\right)$ of order $r$ as,
\begin{eqnarray}
&&H=\left(\mathrm{T}^{(r)}\left(\theta\right)\right)^{-1}_{\theta=0}\left(\partial_\theta
\left(\mathrm{T}^{(r)}\left(\theta\right)\right)\right)_{\theta=0}\\
&&\phantom{H}=\sum_{k=1}^rI\otimes I\otimes\cdots\otimes
\dot{\widehat{R}}_{k,k+1}\left(0\right)\otimes
I\otimes\cdots\otimes I,
\end{eqnarray}
where
\begin{equation}
\dot{\widehat{R}}_{k,k+1}\left(0\right)=\left(\partial_\theta
\widehat{R}_{k,k+1}\left(\theta\right)\right)_{\theta=0}.
\end{equation}
For cyclic boundary condition (and order $r$) one imposes
$k+1=r+1\approx 1$). Thus for $r=2$
\begin{equation}
H=\dot{\widehat{R}}_{12}\left(0\right)+\dot{\widehat{R}}_{21}\left(0\right)=
\dot{\widehat{R}}\left(0\right)+P\dot{\widehat{R}}\left(0\right)P.
\end{equation}
For $\widehat{R}\left(\theta\right)=\frac
1{\sqrt{1+\tanh^2\theta}}\left(I\otimes I+\tanh\theta K\otimes
J\right)$, we obtain
\begin{equation}
\dot{\widehat{R}}\left(0\right)=K\otimes J.
\end{equation}
For the second class, i.e. when
$\widehat{R}\left(\theta\right)=\frac
1{\sqrt{1+\tanh^2\theta}}\left(I\otimes I+\tanh\theta J\otimes
K\right)$, we deduce
\begin{equation}
\dot{\widehat{R}}\left(0\right)=J\otimes K=P\left(K\otimes
J\right)P.
\end{equation}
Thus, in particular, for $r=2$, both classes lead to the same
result
\begin{equation}
H=K\otimes J+J\otimes K=
\sum_{i,j}\left((-1)^{\overline{j}}+(-1)^{\overline{i}}\right)\left((i\overline{i})\otimes
(j\overline{j})\right).
\end{equation}

\section{Potentials for factorizable $S$-matrices}
\setcounter{equation}{0}

Such potentials have been studied in our previous papers
\cite{R8,R9} where basic sources are cited. They are given by
inverse Cayley transforms
\begin{equation}
-\hbox{\verb"i"}\mathrm{V}\left(z\right)=\left(R\left(z\right)-\lambda\left(z\right)I\right)^{-1}\left(R\left(z\right)
+\lambda\left(z\right)I\right),
\end{equation}
where $R\left(z\right)=P\widehat{R}\left(z\right)$ is the
YB-matrix and the parameter $\lambda\left(z\right)$ has been
introduced to guarantee the existence of the inverse \cite{R8,R9}.
To absorb the normalization factor of our unitary
$R\left(z\right)$ we define
\begin{eqnarray}
&&\mu\left(z\right)=\sqrt{1+z^2}\lambda\left(z\right),\\
&&-\hbox{\verb"i"}\mathrm{V}\left(z\right)=\left(\sqrt{1+z^2}R\left(z\right)-\mu\left(z\right)I\right)^{-1}\left(\sqrt{1+z^2}R\left(z\right)
+\mu\left(z\right)I\right)\nonumber\\
&&\phantom{-iV\left(z\right)}=I_{(2n)^2}+2\mu\left(z\right)\left(\sqrt{1+z^2}R\left(z\right)-
\mu\left(z\right)I\right)^{-1}\nonumber\\
&&\phantom{-iV\left(z\right)}\equiv
I_{(2n)^2}+2\mu\left(z\right)X\left(z\right).
\end{eqnarray}
We consider YB-matrices for our class I. For $n=1$,
\begin{equation}
\sqrt{1+z^2}R\left(z\right)=\left|\begin{array}{cccc}
   1 & 0 & 0 & z \\
   0 & z & 1 & 0 \\
   0 & 1 & -z & 0 \\
   -z & 0 & 0 & 1\\
\end{array}\right|.
\end{equation}
Hence (suppressing arguments $z$)
\begin{equation}
X\left(z\right)=\left|\begin{array}{cccc}
   K_1\left(1-\mu\right) & 0 & 0 & -K_1z \\
   0 & K_2\left(z+\mu\right) & K_2 & 0 \\
   0 & K_2 & -K_2\left(z-\mu\right) & 0 \\
   K_1z & 0 & 0 & K_1\left(1-\mu\right)\\
\end{array}\right|,
\end{equation}
where
\begin{equation}
K_1=\frac 1{\left(1-\mu\right)^2+z^2},\qquad K_2=\frac
1{\left(z^2-\mu^2\right)+1}.
\end{equation}
From (8.2) and (8.6) it follows that the inverse $X\left(z\right)$
is well defined for $\lambda\left(z\right)\neq \pm
1,\;\left(\frac{1\pm iz}{1\mp iz}\right)^{1/2}$. With a
$\mathrm{V}$ satisfying (8.1) where $R\left(z\right)$ is a
YB-matrix and the parameter $\mu\left(z\right)$ and hence
$\lambda\left(z\right)$ suitably chosen so that $K_1$, $K_2$ are
well defined the Lagrangian is constructed in terms of the
elements $\mathrm{V}_{(ab,cd)}$ \cite{R10,R11} where
\begin{equation}
\mathrm{V}=\sum_{ab,cd}\mathrm{V}_{(ab,cd)}\left(ab\right)\otimes
\left(cd\right).
\end{equation}
The fermionic Lagrangian is of the form
\begin{equation}
{\cal L}=\int
dx\left(i\overline{\psi}_a\gamma_\nu\partial_\nu\psi_a-g\left(\overline{\psi}_a\gamma_\nu\psi_c\right)
\mathrm{V}_{ab,cd}\left(\overline{\psi}_b\gamma_\nu\psi_d\right)\right).
\end{equation}
The scalar Lagrangian has an interaction term of the form
$\left(\overline{\varphi}_a\overline{\varphi}_c\right)
\mathrm{V}_{ab,cd}\left(\varphi_b\varphi_d\right)$. For $n=2$,
using the same $\left(K_1,K_2\right)$ one obtains the non-zero
elements of $X$ as
\begin{eqnarray}
&&
X\left(1j,1k\right)=K_1\left(\left(1-\mu\right)\delta_{1j}\delta_{1k}-z\delta_{4j}\delta_{4k}\right),\nonumber\\
&&
X\left(4j,4k\right)=K_1\left(\left(1-\mu\right)\delta_{4j}\delta_{4k}+z\delta_{1j}\delta_{1k}\right),\nonumber\\
&&
X\left(1j,4k\right)=K_2\left(\left(\mu+z\right)\delta_{1j}\delta_{4k}+\delta_{4j}\delta_{1k}\right),\nonumber\\
&&
X\left(4j,1k\right)=K_2\left(\left(\mu-z\right)\delta_{4j}\delta_{1k}+\delta_{1j}\delta_{4k}\right),\nonumber\\
&&
X\left(2j,2k\right)=K_1\left(\left(1-\mu\right)\delta_{2j}\delta_{2k}+z\delta_{3j}\delta_{3k}\right),\nonumber\\
&&
X\left(3j,3k\right)=K_1\left(\left(1-\mu\right)\delta_{3j}\delta_{3k}-z\delta_{2j}\delta_{2k}\right),\nonumber\\
&&
X\left(2j,3k\right)=K_2\left(\left(\mu-z\right)\delta_{2j}\delta_{3k}+\delta_{3j}\delta_{2k}\right),\nonumber\\
&&
X\left(3j,2k\right)=K_2\left(\left(\mu+z\right)\delta_{3j}\delta_{2k}+\delta_{2j}\delta_{3k}\right),\nonumber\\
&&
X\left(1j,2k\right)=K_2\left(\mu\delta_{1j}\delta_{2k}+\delta_{2j}\delta_{1k}+z\delta_{3j}\delta_{4k}\right),\nonumber\\
&&
X\left(2j,1k\right)=K_2\left(\delta_{1j}\delta_{2k}+\mu\delta_{2j}\delta_{1k}-z\delta_{4j}\delta_{3k}\right),\nonumber\\
&&
X\left(3j,4k\right)=K_2\left(z\delta_{1j}\delta_{2k}+\mu\delta_{3j}\delta_{4k}+\delta_{4j}\delta_{3k}\right),\nonumber\\
&&
X\left(4j,3k\right)=K_2\left(-z\delta_{2j}\delta_{1k}+\delta_{3j}\delta_{4k}+\mu\delta_{4j}\delta_{3k}\right).
\end{eqnarray}
Finally, defining
\begin{eqnarray}
&&
C_1=\left(\mu\delta_{1j}\delta_{3k}+\delta_{3j}\delta_{1k}-z\delta_{2j}\delta_{4k}\right),\nonumber\\
&&
C_2=\left(\mu\delta_{2j}\delta_{4k}+\delta_{4j}\delta_{2k}+z\delta_{1j}\delta_{3k}\right),\nonumber\\
&&
C_3=\left(\mu\delta_{3j}\delta_{1k}+\delta_{1j}\delta_{3k}-z\delta_{4j}\delta_{2k}\right),\nonumber\\
&&
C_4=\left(\mu\delta_{4j}\delta_{2k}+\delta_{2j}\delta_{4k}+z\delta_{3j}\delta_{1k}\right),\nonumber
\end{eqnarray}
we obtain
\begin{eqnarray}
&&
X\left(1j,3k\right)=K_1K_2\left(K_2^{-1}C_1-2\mu zC_2\right),\nonumber\\
&&
X\left(2j,4k\right)=K_1K_2\left(K_2^{-1}C_2+2\mu zC_1\right),\nonumber\\
&&
X\left(3j,1k\right)=K_1K_2\left(K_2^{-1}C_3-2\mu zC_4\right),\nonumber\\
&& X\left(4j,2k\right)=K_1K_2\left(K_2^{-1}C_4+2\mu zC_3\right).
\end{eqnarray}
Branching through successive scatterings can be two- or three-fold
at each stage. Thus, for example, schematically,
\begin{figure}[ht]
\centerline{\includegraphics[height=9cm]{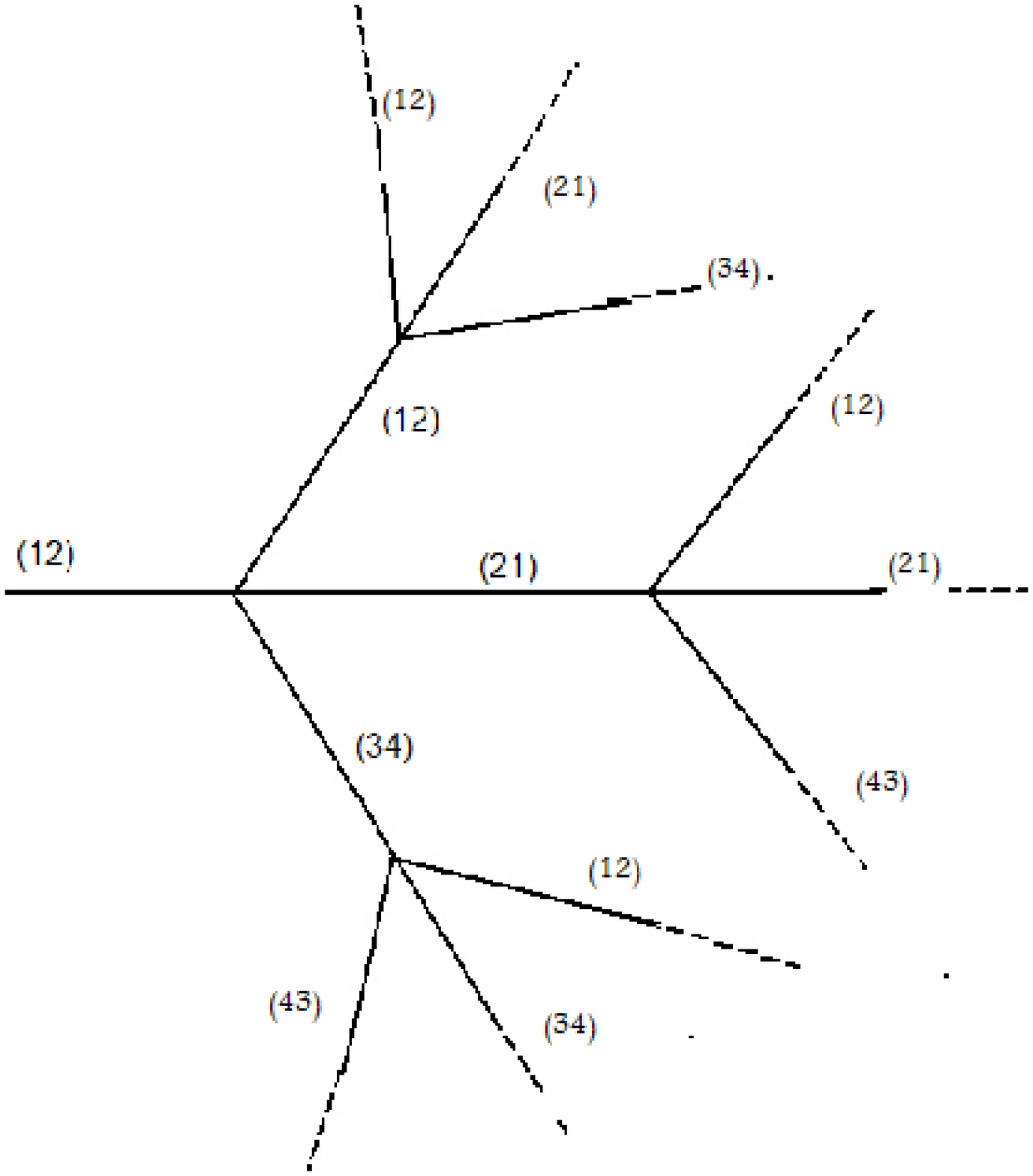}}
\end{figure}

\section{Non-commutative spaces}
\setcounter{equation}{0}

Use of the projectors to construct covariant calculus for
non-commutative spaces has been studied in previous paper
\cite{R1,R7} where basic sources \cite{R12,R13} are cited. We
present below briefly the generalization of the formalism of sec.
2.3 of ref. \cite{R1} to cases $n>1$.

The projectors (3.1), $(2n)^2\times (2n)^2$ matrices,
$P_{\pm}=\frac 12\left(I\otimes I\pm \verb"i"K\otimes J\right)$
being complex, though $\widehat{R}$ is real, special features
arise. Setting, for example, $\mathrm{X}$ being the column vector
of the coordinates,
\begin{equation}
P_{-}\mathrm{X}\otimes \mathrm{X}=0
\end{equation}
one obtains
\begin{equation}
\mathrm{X}_i\mathrm{X}_j=\verb"i"(-1)^{\overline{j}}\mathrm{X}_{\overline{i}}\mathrm{X}_{\overline{j}}.
\end{equation}
Since $\verb"i"(-1)^{\overline{j}}=-\verb"i"(-1)^{j}$ one half of
the constraints repeats, in consistent fashion, the other half.
Thus, for $n=1$
\begin{equation}
X_1X_1=\verb"i"X_{2}X_{2},\qquad X_1X_2=-\verb"i"X_{2}X_{1}
\end{equation}
are sufficient.

Introducing the $(2n)\times (2n)$ projectors $\frac 12\left(I\pm
\verb"i" J\right)$ and defining
\begin{eqnarray}
&&\mathrm{x}_{\pm}=\frac 12\left(I\pm \verb"i" J\right)\mathrm{X},\\
&&2P_-\left(\mathrm{X}\otimes \mathrm{X}\right)=\mathrm{X}\otimes
\mathrm{X}-\left(K\mathrm{X}\right)\otimes\left(\verb"i"J\mathrm{X}\right)\nonumber\\
&&\phantom{2P_-\left(\mathrm{X}\otimes
\mathrm{X}\right)}=\left(\mathrm{x}_+
+\mathrm{x}_-\right)\otimes\left(\mathrm{x}_++\mathrm{x}_-\right)-
\left(\overline{\mathrm{x}}_++\overline{\mathrm{x}}_-\right)\otimes
\left(\mathrm{x}_+-\mathrm{x}_-\right),
\end{eqnarray}
where
\begin{equation}
K\mathrm{X}=\overline{\mathrm{X}}=\left(\overline{\mathrm{x}}_++\overline{\mathrm{x}}_-\right)
\end{equation}
has been implemented. Thus, in terms of $\mathrm{x}_{\pm}$ the
constraints are real.

Define
\begin{equation}
\mathcal{Q}=\nu P_+-I\otimes I,
\end{equation}
where $\nu\neq 0,1$, is otherwise a free parameter to start with.
Then
\begin{equation}
\mathcal{Q}^{-1}=\frac 1{\nu-1}\left(I\otimes I-\nu P_-\right).
\end{equation}
A covariant prescription for the differentials $\mathrm{Z}$ and
the modular structure is
\begin{equation}
\mathcal{Q}\left(\mathrm{Z}\otimes\mathrm{X}\right)=\mathrm{X}\otimes
\mathrm{Z}, \qquad P_+\left(\mathrm{Z}\otimes \mathrm{Z}\right)=0.
\end{equation}
One can define a mobile frame \cite{R7,R14,R15} starting with
\begin{equation}
\mathrm{\Theta}=\sum_i\mathrm{\Theta}_iZ_i,
\end{equation}
where the coefficients $\mathrm{\Theta}_i$ are to be so
constructed that
\begin{equation}
\left[\mathrm{\Theta},\mathrm{X}_i\right]=0.
\end{equation}
A computation leads to
\begin{equation}
\mathrm{X}_i\mathrm{\Theta}_j=\sum_{k,l}\mathrm{\Theta}_k\left(\mathcal{Q}^{-1}P\right)_{kj,il}\mathrm{X}_l,
\end{equation}
where $\mathcal{Q}^{-1}$ is given by (9.8) and the permutation
matrix $P$ by (2.22). We now show how to relate
$\mathcal{Q}^{-1}P$ to $\mathrm{L}^{\pm}$ as in the references
cited above. This involves particular choices of $\nu$. From
(3.3), its easy to see that $\widehat{R}^{\pm 1}=e^{\mp
\verb"i"\pi/4}\left(I\otimes I-\left(1\mp
\verb"i"\right)P_-\right)$. By choosing $\nu=1\mp \verb"i"$, we
can write, using (6.2), that
\begin{equation}
\mathcal{Q}^{-1}P=e^{\mp\mathrm{i}\pi/4}\left(\widehat{R}^{\pm
1}P\right)=e^{\mp\mathrm{i}\pi/4}\mathrm{L}^{\pm}.
\end{equation}
Thus, for example, one can set with $\nu=1+\mathrm{i}$,
\begin{equation}
\mathrm{X}_i\mathrm{\Theta}_j=\sum_{k,l}\mathrm{\Theta}_k\left(e^{\mathrm{i}\pi/4}
\mathrm{L}^{-}\right)_{kj,il}\mathrm{X}_l.
\end{equation}
We will not undertake here any explicit construction of
$\mathrm{\Theta}_i$.

\section{Link invariants (Turaev constructions)}
\setcounter{equation}{0}

We now construct an {\it enhanced system} \cite{R7,R16,R17}
starting with our unitary braid matrices. This implies explicit
construction of a $(2n)\otimes (2n)$ matrix $\mathcal{F}$
satisfying
\begin{eqnarray}
&&\widehat{R}^{\pm 1}\left(\mathcal{F}\otimes
\mathcal{F}\right)=\left(\mathcal{F}\otimes
\mathcal{F}\right)\widehat{R}^{\pm 1},\\
&&\hbox{Tr}_2\left(\widehat{R}^{\pm 1}\left(\mathcal{F}\otimes
\mathcal{F}\right)\right)=a^{\pm 1}b\mathcal{F},\end{eqnarray}
where $(a,b)$ are invertible parameters and
$\hbox{Tr}_2\left(\sum_{i,j,k,l}c_{ij,kl}(ij)\otimes
(kl)\right)=\sum_{i,j}\left(\sum_kc_{ij,kk}\right)(ij)$.

It is sufficient to consider our class I, i.e. $\widehat{R}^{\pm
1}=\frac 1{\sqrt{2}}\left(I\otimes I\pm K\otimes J\right)$. Define
\begin{equation}
\mathcal{F}=\sum_{j=1}^n d_j\left((jj)+(\overline{j}
\overline{j})\right).
\end{equation}
Its follows that
\begin{equation}
J\mathcal{F}=\mathcal{F}J=\sum_{j=1}^n(-1)^{\overline{j}}\left((j\overline{j})-(\overline{j}j)\right)d_j,\qquad
K\mathcal{F}=\mathcal{F}K=\sum_{j=1}^n\left((j\overline{j})+(\overline{j}j)\right)d_j.
\end{equation}
Hence
\begin{equation}
\widehat{R}^{\pm 1}\left(\mathcal{F}\otimes\mathcal{F}\right)
=\frac 1{\sqrt{2}}\left(\mathcal{F}\otimes \mathcal{F}\pm
K\mathcal{F}\otimes J\mathcal{F}\right)=\frac
1{\sqrt{2}}\left(\mathcal{F}\otimes \mathcal{F}\pm
\mathcal{F}K\otimes\mathcal{F}J\right)=\left(\mathcal{F}\otimes
\mathcal{F}\right)\widehat{R}^{\pm 1}.
\end{equation}
Using the results $\hbox{Tr}_2\left(\mathcal{F}\otimes
\mathcal{F}\right)=\mathcal{F}\hbox{Tr}\left(\mathcal{F}\right)=
2\left(\sum_{j=1}^nd_j\right)\mathcal{F}$ and $\hbox{Tr}_2
\left(\left(K\otimes J\right)\left(\mathcal{F}\otimes
\mathcal{F}\right)\right)=\left(K\mathcal{F}\right)\hbox{Tr}
\left(J\mathcal{F}\right)=0$, we obtain
\begin{equation}
\hbox{Tr}_2\left(\widehat{R}^{\pm 1}\mathcal{F}\otimes
\mathcal{F}\right)=\frac 1{\sqrt{2}}\hbox{Tr}_2\left(\mathcal{F}
\otimes\mathcal{F}\pm K\mathcal{F}\otimes J\mathcal{F}\right)
=\frac 1{\sqrt{2}}\hbox{Tr}_2\left(\mathcal{F}\otimes
\mathcal{F}\right)=\sqrt{2}\left(\sum_{j=1}^nd_j\right)\mathcal{F}.
\end{equation}
Thus (10.2) is also satisfied with
\begin{equation}
a=1,\qquad b=\sqrt{2}\left(\sum_{j=1}^nd_j\right).
\end{equation}
Note that for $n=1$, $\mathcal{F}$ degenerates to (a factor times)
the unit matrix $I_2$. From $n=2$ onwards a structure begins to
appear along the diagonal. Thus, for $n=2$,
\begin{equation}
\mathcal{F}=\left|\begin{array}{cccc}
  d_1 & 0 & 0 & 0\\
  0 & d_2 & 0 & 0\\
  0 & 0 & d_2 & 0\\
  0 & 0 & 0 & d_1\\
\end{array}\right|,
\end{equation}
where $d_1=d_2$ is not excluded, but in general $d_1\neq d_2$. The
preceding construction can also be carried through (with the same
$\mathcal{F}$) for the baxterized $\widehat{R}^{\pm
1}\left(z\right)$. {\it For each $n$, $\mathcal{F}$ has $n$ free
parameters.}

With $\left(\mathcal{F},a,b\right)$ thus obtained one can define
(since in our case $a=1$)
\begin{equation}
\wp\left(\beta\right)=b^{(-m+1)}\hbox{Tr}\left(\rho_m\left(\beta\right)\cdot
\mathcal{F}^{\otimes m}\right),
\end{equation}
where $\wp\left(\beta\right)$ is the representation of the braid
$\beta$ associated to $\widehat{R}$ and $\rho_m$ is the
endomorphism of $V^{\otimes m}$. It can be shown (sec. 15 of ref.
\cite{R17}) that this provides an invariant of oriented links,
Markov invariance being assured. For {\it unknot} (no crossing)
\begin{equation}
\wp\left(0\right)=\hbox{Tr}\left(\mathcal{F}\right)=b/\sqrt{2}.
\end{equation}
For our unitary matrices ($I$ denoting $I_{(2n)}$)
\begin{equation}
\widehat{R}+\widehat{R}^{-1}=\sqrt{2}I\otimes I=\sqrt{2}I_{(2n)^2}
\end{equation}
and hence, for all $n$,
\begin{equation}
\widehat{R}^4=-I_{(2n)^2}, \qquad \widehat{R}^8=I_{(2n)^2},
\end{equation}
(one obtains (10.12) most directly by writing (3.3) as
$\widehat{R}=e^{-\mathrm{i}\pi/4}P_++e^{\mathrm{i}\pi/4}P_-$).
Restrictions on the skein relations and the periodicity
(eight-fold) implicit in (10.11) and (10.12) are pointed out in
sec. 4 of ref. \cite{R3}.

\section{Entangled remarks}
\setcounter{equation}{0}

Though, as we explicitly displayed in sec. 4, $\widehat{R}$ block
diagonalized with $\widehat{R}_{(2)}$ (generalization of (4.9) for
$n>2$ being direct) does not satisfy the braid equation, yet one
can write (6.8) as
\begin{eqnarray}
&&\left(V\widehat{R}\left(\theta-\theta'\right)V^{-1}\right)\left(V\mathrm{L}_2\left(\theta\right)V^{-1}\right)
\left(V\mathrm{L}_1\left(\theta'\right)V^{-1}\right)=\left(V\mathrm{L}_2\left(\theta\right)V^{-1}\right)\times\nonumber\\
&&\phantom{xxxxxxxxxxxxxxxx}
\left(V\mathrm{L}_1\left(\theta'\right)V^{-1}\right)\left(V\widehat{R}\left(\theta-\theta'\right)V^{-1}\right).
\end{eqnarray}
Block diagonal ansatz for $\left(V\mathrm{L}_1\left(\theta'\right)
V^{-1}\right)$, $i=1,2$, in terms of the $\mathrm{L}$-functions
for $\mathrm{SO3}$ (studied extensively in \cite{R2}) will
evidently satisfy (11.1). Structures obtained on conjugating back
with $V^{-1}$ might be of interest. But explicit verification is
necessary. One can treat the $\widehat{R}\mathrm{TT}$-relations
analogously. Certain crucial properties of $\mathrm{L}$- and
$\mathrm{T}$-functions have been presented in sec. 6. We hope to
present a more through study elsewhere.

As already pointed out in \cite{R2}, statistical models associated
to unitary $\widehat{R}$ cannot have only non-negative Boltzmann
weights. Negative and complex weights need suitable
interpretations. But the simple structure of the Hamiltonians
governing evolution of the states is worth noting. Our (7.7) is an
example of all $n$. The complex non-commutative spaces associated
to unitary $\widehat{R}$ (sec. 6) deserve further study. So does
the complex unitary braid matrix for odd dimensions (sec. 5). {\it
Topological entanglements} have been presented in \cite{R3} in
terms of link invariants and topological fields. Their possible
relations with quantum entanglements have been emphasized. Here we
have briefly presented Turaev construction of link invariants for
all $n$. Concerning fields we have shown how our $\widehat{R}$ can
be implemented in constructing potentials for factorizable
$S$-matrices. Already, for $n=2$, the structure is considerably
enriched. A canonical formulation for all $n$ would be
interesting.

In \cite{R3} the unitary matrix (1.1) is presented as a common
source of quantum and also of topological entanglements. Acting on
the base space of states
$\left(\left|++\right\rangle,\left|+-\right \rangle,
\left|-+\right\rangle,\left|--\right\rangle\right)$ $\widehat{R}$
of (1.1) generates entangled Bell-states. But $\widehat{R}$ also
has braid property and hence leads to link invariants, links being
viewed as topological entanglements. It also leads to $S$-matrices
where one can permute the successive scattering (sec. 8). Our
$(2n)^2\times (2n)^2$ braid matrices generalize the states
\begin{equation}
\frac 1{\sqrt{2}}\left(\left|++\right\rangle\pm
\left|--\right\rangle\right),\qquad \frac
1{\sqrt{2}}\left(\left|+-\right\rangle\pm
\left|-+\right\rangle\right)
\end{equation}
for (1.1) to
\begin{equation}
\frac 1{\sqrt{2}}\left(I\otimes I+K\otimes
J\right)\left|\mathbb{V}\right\rangle\otimes\left|\mathbb{V}\right\rangle,
\end{equation}
where one may adopt the notation of {\it spin-n} components with
\begin{equation}
\left|\mathbb{V}\right\rangle\equiv\left(\begin{array}{c}
   \left|n\right\rangle\\
  \left|n-1\right\rangle\\
  \vdots\\
  \left|1\right\rangle\\
  \left|-1\right\rangle\\
  \vdots\\
  \left|-n+1\right\rangle\\
  \left|-n\right\rangle\\
\end{array}\right).
\end{equation}
One obtains as direct generalizations of (11.2) the states
\begin{equation}
\frac 1{\sqrt{2}}\left(\left|n-j\right\rangle
\left|n-k\right\rangle\pm \left|-n+j\right\rangle
\left|-n+k\right\rangle\right),
\end{equation}
where $0\leq j,k\leq n-1$. Here the subspaces of
$\left|\mathbb{V}\right\rangle$ can denote any property of the
system under consideration with $2n$ orthogonal states, through
one may use the term {\it spin} for convenience. Now the link
invariants will correspond to more general constructions of sec.
10. The quantum and topological entanglements are both generalized
simultaneously in the sense indicated above. The potential for
factorizable $S$-matrix is now generalized too as indicated in
sec. 8. Note that if instead of using $\widehat{R}$ of (2.7) one
implements the Baxterized $\widehat{R}\left(z\right)$ of (2.12),
one obtains, in place of (11.5) the superpositions
\begin{equation}
\frac 1{\sqrt{1+z^2}}\left(\left|n-j\right\rangle
\left|n-k\right\rangle\pm z\left|-n+j\right\rangle
\left|-n+k\right\rangle\right).
\end{equation}
The matrix $\widehat{R}\left(z\right)$ is still unitary, though
(2.1) is now replaced by (2.19).

Corresponding to odd dimensional and complex (but unitary)
$\widehat{R}\left(\theta\right)$ one has superpositions with
complex coefficients. For simplicity we consider only the
simplest, but typical, case of $n=1$ corresponding to (5.5). Using
a simple, evident, notation
\begin{equation}
\widehat{R}\left(\theta\right)\left(\left(\begin{array}{c}
   \left|+\right\rangle\\
  \left|0\right\rangle\\
  \left|-\right\rangle\\
\end{array}\right)\otimes \left(\begin{array}{c}
   \left|+\right\rangle\\
  \left|0\right\rangle\\
  \left|-\right\rangle\\
\end{array}\right)\right)
\end{equation}
yields superpositions (with $\left(a_{\pm},b_{\pm},c_{\pm}\right)$
of (5.3))
\begin{equation}
a_{\pm}\left|++\right\rangle+a_{\mp}\left|--\right\rangle,\,\,
b_{\pm}\left|+0\right\rangle+b_{\mp}\left|-0\right\rangle,\,\,
a_{\pm}\left|+-\right\rangle+a_{\mp}\left|-+\right\rangle,\,\,
c_{\pm}\left|0+\right\rangle+c_{\mp}\left|+0\right\rangle,\,\,
\left|00\right\rangle
\end{equation}
Note the special status of the central state
$\left|00\right\rangle$. It was explained in ref. \cite{R6} how
the structure of $\widehat{R}\left(\theta\right)$ depends
crucially on the existence of the central element 1. The same
feature singles out $\left|00\right\rangle$ in (11.8). Setting,
say,
\begin{equation}
m_{ij}^{(+)}=-m_{ij}^{(-)}=m_{ij}
\end{equation}
one obtains the simpler superpositions
\begin{equation}
\cos\left(m_{11}\theta\right)\left|\pm\pm\right\rangle+i\sin\left(m_{11}\theta\right)
\left|\mp\mp\right\rangle
\end{equation}
and so on. The content and significance of parameter dependent,
unitary rotations of the base space (11.6), (11.8), (11.10), where
the matrix involved satisfies Baxterized braid equation, deserves
further study.

\vskip 0.5cm

\renewcommand{\theequation}{A.\arabic{equation}}

\paragraph{ADDENDUM:}\setcounter{equation}{0} After having completed this paper we received
the preprint of Yong Zhang and Mo-Lin Ge \cite{Addendum}. In their
construction, in contrast to ours, there are some phases involving
free parameters on the anti-diagonal of the braid matrices. Such
phases provide what the authors designate as {\it generalized}
constructions. But {\it these phases are spurious}. They can be
absorbed by conjugations by matrices of the form $Y\otimes Y$,
where $Y$ is a $\left(2n\right)\times \left(2n\right)$ invertible
matrix for $\widehat{R}$ of dimension $\left(2n\right)^2\times
\left(2n\right)^2$. Such conjugations {\it preserve the braid
property}. The required conjugation is simple. It is sufficient to
illustrate this for the $4\times 4$ case (eq. (29)
\cite{Addendum}). The generalization for $n>1$ will be evident,
recognizing the constraint (28) of ref. \cite{Addendum} as a
crucial feature. Define
\begin{equation}
Y=\left|\begin{array}{cc}
   e^{-i\varphi/4} & 0\\
  0 & e^{i\varphi/4}\\
\end{array}\right|.
\end{equation}
Then
\begin{equation}
Y\otimes Y \left|\begin{array}{cccc}
   0& 0 &0 & e^{i\varphi}\\
0& 0 &1 & 0\\
0& -1 &0 & 0\\
-e^{-i\varphi}& 0 &0 & 0\\
\end{array}\right|Y^{-1}\otimes Y^{-1}=\left|\begin{array}{cccc}
   0& 0 &0 & 1\\
0& 0 &1 & 0\\
0& -1 &0 & 0\\
-1& 0 &0 & 0\\
\end{array}\right|=\left|\begin{array}{cc}
   0& 1\\
-1& 0 \\
\end{array}\right|\otimes \left|\begin{array}{cc}
   0& 1\\
1& 0\\
\end{array}\right|.
\end{equation}
This is $L\otimes K$ (=$J\otimes K$ for $n=1$) of our (4.14). (The
remaining term of $\sqrt{2}\widehat{R}$ is $I\otimes I$ and
invariant.) Removal of the phases makes the tensored structure
($L\otimes K$, for example) evident. Thus one arrives at our
formalism, deriving results in canonical forms with ease and
power.

One can certainly introduce not only phases but a more complicated
parametrization in the tensored structure, trivially, by
substituting
\begin{equation}
J\otimes K\Rightarrow \left(X_{(2n)}\otimes
X_{(2n)}\right)\left(J\otimes K\right)\left(X_{(2n)}^{-1}\otimes
X_{(2n)}^{-1}\right)=\left(X_{(2n)}JX_{(2n)}^{-1}\right)\otimes
\left(X_{(2n)}K X_{(2n)}^{-1}\right)
\end{equation}
and so on, $X_{(2n)}$ being any invertible $\left(2n\right)\times
\left(2n\right)$ matrix. Such parameters are, evidently, to be
removed when present rather than introduced. Unitary
$X_{\left(2n\right)}$ can preserve unitarity, but even so, should
be eliminated.

We would like to add that this is not the first time we {\it
conjugate away} spurious parameters haunting the literature on
braid matrices. A much less simple exercise was provided by the
so-called {\it hybrid deformations} \cite{AACDMplus}. (Another
claim was recently disproved in a note \cite{CDMplus}, using
however a different type of argument.) Being very much conscious
of the possibility of hidden equivalences we scrupulously
displayed the {\it non-equivalence} embodied in our (4.12). This
was, as explained, related to the absence of a tensored structure
($V\neq Y\otimes Y$ for some suitable $Y$).

\vskip 0.5cm

\noindent{\bf Acknowledgments:} {\em One of us (BA) wants to thank
Paul Sorba for precious help. This work is supported by a grant of
CMEP program under number 04MDU615. The work of VKD and SGM was
supported in part by the Bulgarian National Council for scientific
Research, grant F-1205/02 and the European RTN 'Force-universe',
contract MRTN-CT-2004-005104, and by the Alexander Von Humboldt
Foundation in framework of the Clausthal-Leipzig-Sofia
cooperation. We thank Yong Zhang and Mo-Lin Ge for communicating
their paper acknowledged in our Addendum.}


\vskip 0.5cm

\begin{thebibliography}{99}

\bibitem{R1} D. Arnaudon, A. Chakrabarti, V.K. Dobrev and S. Mihov,
{\it Spectral decomposition and Baxterisation of exotic bialgebras
and associated non-commutative geometries}, Int. J. Mod. Phys.
{\bf A18} (2003) 4201. \texttt{math.QA/0209321}.

\bibitem{R2} D. Arnaudon, A. Chakrabarti, V.K. Dobrev and S. Mihov,
{\it Exotic bialgebra $SO3$: Representations, Baxterization and
Applications}, Ann. H. Poincaré. {\bf 7} (2006) 1351.
\texttt{math.QA/0601708}.

\bibitem{R3} L.H. Kauffman and S.J. Lomonaco Jr., {\it Braiding operators are universal quantum gates},
New. J. Phys. {\bf 6} (2004) 34. \texttt{quant-ph/0401094}.

\bibitem{R4} J. Franko, E.C. Rowell and Z. Wang, {\it Extaspecial
2-groups and images of braid group representations}, Knot Theory
Ramifications, {\bf 15} $N^o$. 4 (2006) 413.
\texttt{math.RT/0503435}.

\bibitem{R5} Y. Zhang, N. Jing, M.L. Ge, {\it Quantum algebras associated with Bell states},
\texttt{math-ph/0610036}.

\bibitem{SB} J.K. Slingerland and F.A. Bais, {\it Quantum groups and non-Abelian braiding in quantum Hall
systems}, Nucl. Phys. B 612 (2001) 229-290.

\bibitem{ZKG} Y. Zhang, L.H. Kauffman and M.L. Ge, {\it Universal quantum gate,
Yang-baxterization and hamiltonian}, Int. J. Quant. Inf. 4 (2005)
669.

\bibitem{R6} A. Chakrabarti, {\it A nested sequence of projectors and corresponding braid
matrices $\hat R(\theta)$: (1) Odd dimensions}, Jour. Math. Phys.
{\bf 46} (2005) 063508. \texttt{math.QA/0401207}.

\bibitem{R7} A. Chakrabarti, {\it Aspects of a new class of braid matrices: Roots of unity and hyperelliptic $q$
for triangularity, $L$-algebra, link-invariants, non-commutatives
spaces}, Jour. Math. Phys. {\bf 46} (2005) 063509.
\texttt{math.QA/0412549}. .

\bibitem{R8} B. Abdesselam and A. Chakrabarti, {\it A nested sequence of projectors: (2) Multiparameter
multistate statistical models, Hamiltonians, $S$-matrices}, Jour.
Math. Phys. {\bf 47}(2006) 053508. \texttt{math.QA/0601584}

\bibitem{R9} B. Abdesselam and A. Chakrabarti, {\it New class o$^N$ of statistical models: Transfer matrix
eigenstates, chain Hamiltonians, factorizable S matrix}, Jour.
Math. Phys. {\bf 47} (2006) 123301. \texttt{math.QA/0607379}.

\bibitem{ACDM} B. Abdesselam, A. Chakrabarti, V.K. Dobrev and S.G. Mihov, (in
preparation).

\bibitem{ACDM1} D. Arnaudon, A. Chakrabarti, V.K. Dobrev and S. Mihov, {\it Duality and representations for new
exotic bialgebras}, J. Math. Phys. {\bf 43} (2002) 6238.
\texttt{math.QA/0206053}.

\bibitem{R10} P.P. Kulish and E.K. Sklyanin, {\it Integrable Quantum Field theories},
Lecture Notes in Physics, Springer, New York, (1982) p. 61.

\bibitem{R11} H.J. De Vega, {\it Yang-Baxter algebras, integrable theories and quantum groups},
Int. Jour. Mod. Phys. A vol. {\bf 4} (1989) 2371.

\bibitem{R12} J. Wess and B. Zumino, Nucl. Phys. B (Proc. Suppl.) {\bf 18B}
(1990) 302.

\bibitem{R13} L. Hlavaty, {\it Yang-Baxter matrices and differential calculi on quantum hyperplanes},
J. Phys. A: Math. Gen. {\bf 25} (1992) 485.

\bibitem{R14} J. Madore, {\it A introduction to noncommutative geometry}, (Cambridge University Press, 1999).

\bibitem{R15} B.L. Cerchiai, G. Fiore and J. Madore, {\it Geometrical Tools for Quantum Euclidean Spaces},
Commun. Math. Phys. {\bf 217} (2001) 521.
\texttt{math.QA/0002007}.

\bibitem{R16} V.G. Tuarev, Invent. Math. {\bf 92} (1988), 527.

\bibitem{R17} V. Chari and A. Pressley, {\it Quantum groups}, Cambridge University Press,
Cambridge, 1994 (sec. 15.1-15.2).

\bibitem{Addendum} Y. Zhang and M.L. Ge, {\it GHZ States, Almost-Complex Structure and
Yang--Baxter Equation (I)}, \texttt{quant-ph/0701244}.

\bibitem{AACDMplus} B.L. Aneva, D. Arnaudon, A. Chakrabarti, V.K. Dobrev and S.
Mihov, {\it On combined standard-nonstandard or hybrid
$\left(q,h\right)$-deformations}, J. Math. Phys. {\bf 42} (2001)
1236.

\bibitem{CDMplus} A. Chakrabarti, V.K. Dobrev and S.
Mihov, {\it On a ''New'' Deformation of GL(2)},
\texttt{math.QA/0701103}.


\end{thebibliography}
\end{document}